\newcount\ite\ite=1\def\0{\global\ite=1\1} 
\def\1{\item{\rm(\romannumeral\the\ite)}\advance\ite1\quad}

\documentclass[12pt,twoside]{article}
\usepackage{amssymb}

\font\teneufm=eufm10 scaled \magstep1
\font\seveneufm=eufm7 scaled \magstep1
\font\fiveeufm=eufm5  scaled \magstep1
\newfam\eufmfam

\textfont\eufmfam=\teneufm
\scriptfont\eufmfam=\seveneufm
\scriptscriptfont\eufmfam=\fiveeufm

\newfam\msbfam
\font\tenmsb=msbm10 scaled \magstep1  \textfont\msbfam=\tenmsb
\font\sevenmsb=msbm7 scaled \magstep1 \scriptfont\msbfam=\sevenmsb
\font\fivemsb=msbm5 scaled \magstep1  \scriptscriptfont\msbfam=\fivemsb

\def\dd#1{\raise1.5pt\hbox{$\,\partial\!$}/\raise-2.5pt\hbox{$\!\partial#1\,$}}

\def\5#1{{\mathcal #1}}

\def\RR{{\mathbb R}}
\def\CC{{\mathbb C}}

\def\ZZ{{\mathbb Z}}

\def\PP{{\mathbb P}}

\def\ra{\rightarrow}

\def\SO{\mathop{\rm SO}\nolimits}

\def\Im{\mathop{\rm Im}\nolimits}
\def\Re{\mathop{\rm Re}\nolimits}

 \def\HollowBoxx #1#2#3{{\dimen0=#1 \advance\dimen0 by -#2
       \dimen1=#1 \advance\dimen1 by #3
        \vrule height 0pt depth #3 width #2
       \hskip -#3
       \vrule height #1 depth #3 width #3}}
 \def\LeftContraction{\mathord{\kern1.45pt \HollowBoxx{6pt}{3.5pt}{.4pt}}\,}

 \def\HollowBox #1#2#3{{\dimen0=#1 \advance\dimen0 by -#3
       \dimen1=#1 \advance\dimen1 by #3
        \vrule height #1 depth #3 width #3
        \vrule height 0pt depth #3 width #2
        \hskip -#3}}
 \def\RightContraction{\mathord{\, \HollowBox{6pt}{3.1pt}{.4pt}} \kern1.6pt}

\def\qed{{\hfill $\Box$}}
\newtheorem{theorem}{THEOREM}[section]
\newtheorem{corollary}[theorem]{Corollary}

\newtheorem{remark}[theorem]{Remark}
\newtheorem{proposition}[theorem]{Proposition}
\newtheorem{conjecture}[theorem]{Conjecture}

  \renewenvironment{thebibliography}[1]{%
    \begin{oldthebibliography}{#1}%
      \setlength{\parskip}{0.1ex}%
      \setlength{\itemsep}{0.5ex}%
  }%
  {%
    \end{oldthebibliography}%
}

\makeatletter
\def\blfootnote{\xdef\@thefnmark{}\@footnotetext}
\makeatother

\begin{document}

\begin{center}


\end{center}

\begin{center}
{\Large \bf On the Classification of Homogeneous\\
\vspace{0.3cm}
Hypersurfaces in Complex Space}
\medskip\medskip\\
\normalsize A. V. Isaev 
\end{center}

\begin{quotation} \small \sl\noindent We discuss a family $M_t^n$, with $n\ge 2$, $t>1$, of real hypersurfaces in a complex affine $n$-dimensional quadric arising in connection with the classification of homogeneous compact simply-connected real-analytic hypersurfaces in\, $\CC^n$ due to Morimoto and Nagano. To finalize their classification, one needs to resolve the problem of the embeddability of $M_t^n$ in\, $\CC^n$ for $n=3,7$. We show that $M_t^7$ is not embeddable in\, $\CC^7$ for every $t$ and that $M_t^3$ is embeddable in\, $\CC^3$ for all $1<t<1+10^{-6}$. As a consequence of our analysis of a map constructed by Ahern and Rudin, we also conjecture that the embeddability of $M_t^3$ takes place for all\, $1<t<\sqrt{(2+\sqrt{2})/3}$. 
\end{quotation}

\thispagestyle{empty}

\pagestyle{myheadings}
\markboth{A. V. Isaev}{Classification of Homogeneous Hypersurfaces}

\setcounter{section}{0}

\section{Introduction}\label{intro}
\setcounter{equation}{0}

For $n\ge 2$, consider the $n$-dimensional affine quadric in $\CC^{n+1}$:
\begin{equation}
Q^n:=\{(z_1,\dots,z_{n+1})\in\CC^{n+1}:z_1^2+\dots+z_{n+1}^2=1\}.\label{quadric}
\end{equation}
The group $\SO(n+1,\RR)$ acts on $Q^n$, and the orbits of this action are the sphere $S^n=Q^n\cap\RR^{n+1}$ as well as the compact strongly pseudoconvex hypersurfaces
\begin{equation}
M_t^n:=\{(z_1,\dots,z_{n+1})\in\CC^{n+1}: |z_1|^2+\dots+|z_{n+1}|^2=t\}\cap Q^n,\,\, t>1.\label{mtn}
\end{equation}
These hypersurfaces play an important role in the classical paper \cite{MN}, where the authors set out to determine all compact simply-connected real-analytic hypersurfaces in $\CC^n$ homogeneous under an action of a Lie group by CR-transformations. They showed that every such hypersurface is CR-equivalent to either the sphere $S^{2n-1}$ or, for $n=3,7$, to the manifold $M_t^n$ for some $t$. However, the question of the existence of a real-analytic CR-embedding of $M_t^n$ in $\CC^n$ for $n=3,7$ was not clarified, thus the classification in these two dimensions was not fully completed.

In this paper, we first discuss the family $M_t^n$ for all $n\ge 2$, $t>1$ (see Section \ref{sect1}). We observe that a necessary condition for the existence of a real-analytic CR-embedding of $M_t^n$ in $\CC^n$ is the embeddability of the sphere $S^n$ in $\CC^n$ as a totally real submanifold (see Corollary \ref{cor}). The problem of the existence of a totally real embedding of $S^n$ in $\CC^n$ was considered by Gromov (see \cite{Gr1} and p. 193 in \cite{Gr2}), Stout-Zame (see \cite{SZ}), Ahern-Rudin (see \cite{AR}), Forstneri\v c (see \cite{F1}--\cite{F3}). In particular, it has turned out that $S^n$ admits a totally real embedding in $\CC^n$ only for $n=3$, hence Corollary \ref{cor} implies that $M_t^n$ cannot be real-analytically CR-embedded in $\CC^n$ for all $t$ if $n\ne 3$. For $n\ne 3,7$ the non-embeddability of $M_t^n$ was established in \cite{MN} by a different method, whereas for $n=7$ it appears to be a new result (see Corollary \ref{cor2}). Furthermore, since $S^3$ is a totally real submanifold of $Q^3$, any real-analytic totally real embedding of $S^3$ into $\CC^3$ extends to a biholomorphic map defined in a neighborhood of $S^3$ in $Q^3$. Owing to the fact that $M_t^3$ accumulate to $S^3$ as $t\ra 1$, this observation implies that $M_t^3$ admits a real-analytic CR-embedding in $\CC^3$ for all $t$ sufficiently close to 1. Thus, the classification of homogeneous compact simply-connected real-analytic hypersurfaces in complex dimension 3 is special as it includes manifolds other than the sphere $S^5$.

Next, in \cite{AR} an explicit polynomial totally real embedding $f$ of $S^3$ into $\CC^3$ was constructed, and in Section \ref{sect2} we study the holomorphic continuation $F:\CC^4\ra\CC^3$ of $f$ in order to determine the interval of parameter values $1<t<t_0$ for which the extended map yields an embedding of $M_t^3$ into $\CC^3$. It turns out that the interval is rather small, namely $t_0\le\sqrt{(2+\sqrt{2})/3}\approx 1.07$ (see Proposition \ref{embed}). A lower bound for $t_0$ is given in Theorem \ref{main}, where we show that $t_0\ge 1+10^{-6}$. Further, our discussion of the fibers of the map $F$ leads to the conjecture asserting that one can in fact take $t_0=\sqrt{(2+\sqrt{2})/3}$ (see Conjecture \ref{conj}). For all other values of $t$ the problem of the embeddability of $M_t^3$ remains completely open. 

{\bf Acknowledgement.} We would like to thank F. Forstneri\v c, X. Huang,\linebreak A. Huckleberry and S. Nemirovskii for useful discussions. This work is supported by the Australian Research Council.

\section{The family $M_t^n$}\label{sect1}
\setcounter{equation}{0}

We start by discussing the family $M_t^n$ for all $n\ge 2$, $t>1$. First of all, we note that the hypersurfaces in the family are all pairwise CR-nonequivalent (see Example 13.9 in \cite{KZ}). Next, computing the isotropy subgroups of the $\SO(n+1,\RR)$-action on $M_t^n$, one observes that $M_t^n$ is diffeomorphic to $\SO(n+1,\RR)/\SO(n-1,\RR)$. On the other hand, from (\ref{quadric}), (\ref{mtn}) we see that
$$
M_t^n=\left\{x+iy\in\CC^{n+1}: ||x||=\sqrt{\frac{t+1}{2}},\, ||y||=\sqrt{\frac{t-1}{2}},\, (x,y)=0\right\},
$$
where $x,y\in\RR^{n+1}$, hence $M_t^n$ is diffeomorphic to the tangent sphere bundle over $S^n$. It then follows that $\pi_1(M_t^n)=0$ if $n\ge 3$, $\pi_1(M_t^2)\simeq\ZZ_2$, and that $M_t^n$ accumulate to $S^n$ as $t\ra 1$. Note also that $M_t^2$ is a double cover of the following well-known homogeneous hypersurface in $\CC\PP^2$ discovered by \'E. Cartan (see \cite{C}): 
$$
\left\{(\zeta_1:\zeta_2:\zeta_3)\in\CC\PP^2:|\zeta_1|^2+|\zeta_2|^2+|\zeta_3|^2=t|\zeta_1^2+\zeta_2^2+\zeta_3^2|\right\}.
$$

For our arguments below we will utilize the strongly pseudoconvex Stein domains in $Q^n$ bounded by $M_t^n$:
$$
D_t^n:=\{(z_1,\dots,z_{n+1})\in\CC^{n+1}: |z_1|^2+\dots+|z_{n+1}|^2<t\}\cap Q^n,\,\, t>1.
$$
The sphere $S^n$ lies in $D_t^n$ for all $t$ and is a strong deformation retract of $D_t^n$. Indeed, the following map from $D_t^n\times[0,1]$ to $D_t^n$ is a strong deformation retraction of $D_t^n$ to $S^n$:
$$ 
(x+iy, s)\mapsto \frac{\sqrt{1+s^2||y||^2}}{||x||}x+isy,\quad 0\le s\le 1.
$$
In particular, $D_t^n$ is simply-connected but not contractible.

\begin{remark}\label{rem1}\rm The domains $D_t^n$ illustrate the relationship between the fundamental group of a smoothly bounded Stein domain and that of its (connected) boundary, which is important, for example, for the uniformization problem. Namely, if $D$ is a smoothly bounded Stein domain, the fundamental groups $\pi_1(D)$ and $\pi_1(\partial D)$ are isomorphic in complex dimensions $\ge 3$, whereas in dimension 2 there exists a surjective homomorphism $\pi_1(\partial D)\ra\pi_1(D)$ and the fundamental group of $\partial D$ can be larger than that of $D$ (see \cite{NS} for a detailed discussion of these facts). Indeed, as we observed above, $\pi_1(D_t^n)=0$ for $n\ge 2$, $\pi_1(M_t^n)=0$ for $n\ge 3$, but $\pi_1(M_t^2)\simeq\ZZ_2$. Examples of simply connected (in fact contractible) domains with non-simply connected boundaries exist even in the class of strongly pseudoconvex domains in $\CC^2$ although they are much harder to construct (see \cite{Go}).
\end{remark}

We will now turn to the problem of the CR-embeddability of $M_t^n$ in $\CC^n$ and prove the following proposition.

\begin{proposition}\label{induced}\sl Any real-analytic CR-embedding of $M_t^n$ into $\CC^n$ extends to a biholomorphic mapping of $D_t^n$ onto a domain in $\CC^n$.
\end{proposition}

\noindent {\bf Proof:} Let $\varphi: M_t^n\ra\CC^n$ be a real-analytic CR-embedding. Then $\varphi$ extends to a biholomorphic map between a neighborhood of $M_t^n$ in $Q^n$ and a neighborhood of the real-analytic strongly pseudoconvex hypersurface $M:=\varphi(M_t^n)$ in $\CC^n$. Further, $\varphi$ extends to a holomorphic map from $D_t^n$ to $\CC^n$ (which follows, for instance, from results of \cite{KR}), and we denote the extension of $\varphi$ to a neighborhood of $\overline{D_t^n}$ also by $\varphi$. 

Next, let $D$ be the strongly pseudoconvex domain in $\CC^n$ bounded by $M$ and $\psi: M\ra M_t^n$ the inverse of $\varphi$ on $M_t^n$. As before, the map $\psi$ extends to a biholomorphic map between a neighborhood of $M$ in $\CC^n$ and a neighborhood of $M_t^n$ in $Q^n$. Furthermore, $\psi$ extends to a holomorphic map from $D$ to $\CC^{n+1}$. We denote the resulting extension of $\psi$ to a neighborhood of $\overline{D}$ also by $\psi$. Clearly, the range of $\psi$ lies in $Q^n$.

It now follows that $\varphi(D_t^n)=D$, $\psi(D)=D_t^n$ and $\varphi$, $\psi$ are the inverses of each other. In particular, $\varphi$ maps $D_t^n$ biholomorphically onto $D$, which completes the proof. \qed
\vspace{0.3cm}

\noindent As an immediate consequence of Proposition \ref{induced} we obtain the following result.

\begin{corollary}\label{cor} \sl If for some $t>1$ the manifold $M_t^n$ is real-analytically CR-embed\-dable in $\CC^n$, then $S^n$ admits a real-analytic totally real embedding\linebreak in $\CC^n$.
\end{corollary}

It is not hard to see that for $n\ne 3,7$ the sphere $S^n$ does not admit a smooth totally real embedding in $\CC^n$ (see \cite{Gr2}, \cite{SZ}). Indeed, multiplication by $i$ establishes an isomorphism between the tangent and the normal bundles of any smooth totally real $n$-dimensional submanifold of $\CC^n$. On the other hand, the normal bundle to $S^n$ induced by any smooth embedding in $\RR^{2n}$ is trivial (see Theorem 8.2 in \cite{K}). Therefore, if $S^n$ admits a smooth totally real embedding in $\CC^n$, it is parallelizable, which is impossible unless $n=3$ or $n=7$. Corollary \ref{cor} then yields that $M_t^n$ cannot be real-analytically CR-embedded in $\CC^n$ for $n\ne 3,7$. This last statement was obtained in \cite{MN} by utilizing the facts that $M_t^n$ is diffeomorphic to the sphere bundle over $S^n$ and that a real-analytic embedding of $M_t^n$ into $\CC^n$ induces a smooth embedding of the tangent bundle of $S^n$ into $\RR^{2n}$, which again leads to the parallelizability of $S^n$ (at least for $n\ge 3$).

Further, the non-existence of a smooth totally real embedding of $S^7$ in $\CC^7$ was first obtained in \cite{SZ} by an argument relying on a result of \cite{S}, which states that any two smooth embeddings of $S^n$ in $\RR^{2n}$ are regularly homotopic. Corollary \ref{cor} then yields: 

\begin{corollary}\label{cor2}\sl No manifold in the family $M_t^7$ admits a real-analytic CR-embedding in\, $\CC^7$, hence every homogeneous compact simply-connected real-analytic hypersurface in $\CC^7$ is CR-equivalent to $S^{13}$. 
\end{corollary}

In contrast, it turns out that $S^3$ can be embedded into $\CC^3$ by a real-analytic CR-map. The first proof of this fact was given in \cite{AR}, where an explicit example of an embedding was constructed. Since $S^3$ is a totally real submanifold of $Q^3$, any real-analytic totally real embedding of $S^3$ into $\CC^3$ extends to a biholomorphic map defined in a neighborhood of $S^3$ in $Q^3$. The manifolds $M_t^3$ accumulate to $S^3$ as $t\ra 1$, hence $M_t^3$ admits a real-analytic CR-embedding in $\CC^3$ for all $t$ sufficiently close to 1. This shows that, interestingly, the classification of homogeneous compact simply-connected real-analytic hypersurfaces in $\CC^3$ includes manifolds other than $S^5$. The embedding found in \cite{AR} is a polynomial map on $\RR^4\subset\CC^4$, hence it has a (unique) holomorphic continuation to all of $\CC^4$. We will study this extended map in the next section.

\begin{remark}\label{rem2}\rm Observe that every hypersurface $M_t^n$ is non-spherical. Indeed, otherwise by results of \cite{NS} the universal cover of the domain $D_t^n$ would be biholomorphic to the unit ball $B^n\subset\CC^n$. Since $D_t^n$ is simply-connected, this would imply that $D_t^n$ is biholomorphic to $B^n$, which is impossible since $D_t^n$ is not contractible. Now, the non-sphericity and homogeneity of the hypersurface $M_t^n$ yield that it has no umbilic points. Therefore, every manifold $M_t^3$ embeddable in $\CC^3$ provides an example of a compact strongly pseudoconvex simply-connected hypersurface in $\CC^3$ without umbilic points. Such hypersurfaces have been known before, but the arguments required to obtain non-umbilicity for them are much more involved than the one given above. For example, the proof in \cite{W} of the fact that every generic ellipsoid in $\CC^n$ for $n\ge 3$ has no umbilic points relies on the Chern-Moser theory (note for comparison that every ellipsoid in $\CC^2$ has at least four umbilic points -- see \cite{HJ}).  
\end{remark}

\section{The holomorphic continuation of the\\ Ahern-Rudin map}\label{sect2}
\setcounter{equation}{0}

In this section we study the holomorphic continuation of the totally real embedding of $S^3$ into $\CC^3$ constructed in \cite{AR}. Let $(z,w)$ be coordinates in $\CC^2$ and let $S^3$ be realized in the standard way as the subset of $\CC^2$ given by
$$
S^3=\{(z,w)\in\CC^2: |z|^2+|w|^2=1\}.
$$
The Ahern-Rudin map is defined on all of $\CC^2$ as follows:
\begin{equation}
f:\CC^2\ra\CC^3,\quad f(z,w):=(z,w,w\overline{z}\overline{w}^2+iz\overline{z}^2\overline{w}).\label{mapf}
\end{equation}
Now, consider $\CC^4$ with coordinates $z_1,z_2,z_3,z_4$ and embed $\CC^2$ into $\CC^4$ as the totally real subspace $\RR^4$:
$$
(z,w)\mapsto(\Re z,\Im z,\Re w,\Im w).
$$
Clearly, the push-forward of the polynomial map $f$ extends from $\RR^4$ to a holomorphic map $F:\CC^4\ra\CC^3$ by the formula
$$
\begin{array}{l}
\displaystyle F(z_1,z_2,z_3,z_4):=\Bigl(z_1+iz_2,z_3+iz_4,\\
\vspace{-0.5cm}\\
\hspace{2.1cm}\displaystyle(z_3+iz_4)(z_1-iz_2)(z_3-iz_4)^2+i(z_1+iz_2)(z_1-iz_2)^2(z_3-iz_4)\Bigr).
\end{array}
$$
It will be convenient for us to argue in the coordinates
\begin{equation}
w_1:=z_1+iz_2,\,\,w_2:=z_1-iz_2,\,\,w_3:=z_3+iz_4,\,\,w_4:=z_3-iz_4.\label{coordw}
\end{equation}
In these coordinates 
the quadric $Q^3$ takes the form
\begin{equation}
\left\{(w_1,w_2,w_3,w_4)\in\CC^4: w_1w_2+w_3w_4=1\right\},\label{qnew}
\end{equation}
the sphere $S^3\subset Q^3$ the form
\begin{equation}
\left\{(w_1,w_2,w_3,w_4)\in\CC^4: w_2=\bar w_1,\,w_4=\bar w_3\right\}\cap Q^3,\label{sphereq}
\end{equation}
the hypersurface $M_t^3\subset Q^3$ the form
\begin{equation}
\left\{(w_1,w_2,w_3,w_4)\in\CC^4: |w_1|^2+|w_2|^2+|w_3|^2+|w_4|^2=2t\right\}\cap Q^3,\label{formm7}
\end{equation}
and the map $F$ the form
\begin{equation}
\displaystyle (w_1,w_2,w_3,w_4)\mapsto\Bigl(w_1,w_3,w_2w_3w_4^2+iw_1w_2^2w_4\Bigr).\label{map}
\end{equation}

We will study the map $F$ in order to obtain some evidence regarding the values of $t$ for which the manifold $M_t^3$ is embeddable in $\CC^3$. Clearly, $F$ defines an embedding of $M_t^3$ if its restriction $\tilde F:=F|_{Q^3}$ is non-degenerate and injective on $M_t^3$, therefore it is important to investigate the non-degeneracy and injectivity properties of $\tilde F$. First of all, observe that $F$ has maximal rank at every point of $Q^3$ since its Jacobian matrix is
$$
\left(\begin{array}{cccc}
1 & 0 & 0 & 0\\
0 & 0 & 1 & 0\\
iw_2^2w_4 & w_3w_4^2+2iw_1w_2w_4 & w_2w_4^2 & 2w_2w_3w_4+iw_1w_2^2
\end{array}
\right),
$$
and the entries $w_3w_4^2+2iw_1w_2w_4$, $ 2w_2w_3w_4+iw_1w_2^2$ cannot simultaneously vanish on $Q^3$. However, as the following proposition shows, most manifolds $M_t^3$ contain points at which the restricted map $\tilde F$ degenerates.

\begin{proposition}\label{embed}\sl The map $\tilde F$ degenerates at some point of $M_t^3$ if and only if $t\ge\sqrt{(2+\sqrt{2})/3}$.  
\end{proposition}

\noindent{\bf Proof:} Observe that $|w_1|+|w_3|>0$ on $Q^3$. For $w_1\ne 0$ we choose $w_1,w_3,w_4$ as local coordinates on $Q^3$ and write the third component of $\tilde F$ as
$$
\varphi:=\frac{1-w_3w_4}{w_1}(iw_4+(1-i)w_3w_4^2)
$$
(see (\ref{qnew}), (\ref{map})). Then the Jacobian $J_{\tilde F}$ of $\tilde F$ is equal to
\begin{equation}
\frac{\partial\varphi}{\partial w_4}=\frac{(3i-3)w_3^2w_4^2+(2-4i)w_3w_4+i}{w_1},\label{phieq}
\end{equation}
hence it vanishes if and only if
$$
w_3w_4=\frac{3\pm\sqrt{2}-i}{6}.
$$
At such points we have
\begin{equation}
|w_1|^2+|w_2|^2+|w_3|^2+|w_4|^2=|w_1|^2+\frac{2\mp\sqrt{2}}{6|w_1|^2}+|w_3|^2+\frac{2\pm\sqrt{2}}{6|w_3|^2}.\label{est1}
\end{equation}

Analogously, if $w_3\ne 0$, we choose $w_1,w_2,w_3$ as local coordinates on $Q^3$ and write the third component of $\tilde F$ as
$$
\psi:=\frac{1-w_1w_2}{w_3}(w_2+(i-1)w_1w_2^2)
$$
(see (\ref{qnew}), (\ref{map})). Then
$$
J_{\tilde F}=-\frac{\partial\psi}{\partial w_2}=-\frac{(3-3i)w_1^2w_2^2+(2i-4)w_1w_2+1}{w_3},
$$
which vanishes if and only if
$$
w_1w_2=\frac{3\pm\sqrt{2}+i}{6}.
$$
Hence for all points of degeneracy of $\tilde F$ we have $w_1\ne 0$, $w_3\ne 0$, and therefore such points are described as the zeroes of ${\partial\varphi}/{\partial w_4}$ or, equivalently, as the zeroes of ${\partial\psi}/{\partial w_2}$.

Now, investigating the behavior of the function
$$
x+\frac{2+\sqrt{2}}{6x}+y+\frac{2-\sqrt{2}}{6y}
$$
for $x,y>0$, one can deduce from (\ref{est1}) that $J_{\tilde F}$ vanishes at a point of $M_t^3$ if and only if $t\ge\sqrt{(2+\sqrt{2})/3}$ as claimed.\qed
\vspace{0.3cm}

As we noted above, any real-analytic totally real embedding of $S^3$ in $\CC^3$ yields an embedding of $M_t^3$ for all $t$ sufficiently close to 1. Define
$$
t_0:=\sup\{t: \hbox{$\tilde F|_{M_s^3}$ is an embedding for all $1<s\le t$}\}.
$$
Proposition \ref{embed} implies $t_0\le\sqrt{(2+\sqrt{2})/3}$. We will now give a lower bound for $t_0$.

\newpage

\begin{theorem}\label{main}\sl One has $t_0\ge 1+10^{-6}$.
\end{theorem}

\noindent{\bf Proof:} For $0<\varepsilon<1$ define
$$
U_{\varepsilon}:=\left\{(w_1,w_2,w_3,w_4)\in\CC^4: |w_2-\bar w_1|<\varepsilon,\,|w_4-\bar w_3|<\varepsilon\right\}\cap Q^3.
$$
Clearly, $U_{\varepsilon}$ is a neighborhood of $S^3$ in $Q^3$ (see (\ref{sphereq})). We will find $\varepsilon$ such that $\tilde F$ is biholomorphic on $U_{\varepsilon}$. Writing any point in $Q^3$ as
$$
W=(w_1,\bar w_1+\mu,w_3,\bar w_3+\eta),
$$
from (\ref{qnew}) we observe
\begin{equation}
|w_1|^2+|w_3|^2+\mu w_1+\eta w_3=1,\label{formmm18}
\end{equation}  
which implies
\begin{equation}
|w_1|^2+|\bar w_1+\mu|^2+|w_3|^3+|\bar w_3+\eta|^2=2+|\mu|^2+|\eta|^2.\label{formm8}
\end{equation}
Clearly, $W$ lies in $U_{\varepsilon}$ if and only if $|\mu|<\varepsilon$, $|\eta|<\varepsilon$. It then follows from (\ref{formm7}), (\ref{formm8}) that $M_t^3\subset U_{\varepsilon}$ for all $1<t<1+\varepsilon^2/2$. Below we will see that one can take $\varepsilon=\sqrt{2}\cdot 10^{-3}$, which will then imply the theorem.

In order to choose $\varepsilon$, we study the fibers of the map $\tilde F$. Let two points $W=(w_1,w_2,w_3,w_4)$ and $\hat W=(\hat w_1,\hat w_2,\hat w_3,\hat w_4)$ lie in $Q^3$ and assume that $\tilde F(W)=\tilde F(\hat W)$. Then from (\ref{map}) one immediately has $\hat w_1=w_1$, $\hat w_3=w_3$ and
\begin{equation}
\hat w_2w_3\hat w_4^{2}+iw_1\hat w_2^{2}\hat w_4=w_2w_3w_4^2+iw_1w_2^2w_4.\label{mainids}
\end{equation}    
If $w_1=0$ or $w_3=0$, then (\ref{mainids}) implies $\hat W=W$, so we suppose from now on that $w_1\ne 0$ and $w_3\ne 0$. Then, using (\ref{qnew}) we substitute
\begin{equation}
w_2=\frac{1-w_3w_4}{w_1},\quad \hat w_2=\frac{1-w_3\hat w_4}{w_1}\label{expr234}
\end{equation}
into (\ref{mainids}) and simplifying the resulting expression obtain
\begin{equation}
\begin{array}{l}
(\hat w_4-w_4)\Biggl[(i-1)w_3^2\,\hat w_4^2+ \Bigl((1-2i)w_3+(i-1)w_3^2w_4\Bigr)\hat w_4+\\
\vspace{-0.3cm}\\
\hspace{4.5cm}\Bigl(i+(1-2i)w_3w_4+(i-1)w_3^2w_4^2\Bigr)\Biggr]=0.\label{mainids4}
\end{array}
\end{equation}
We treat identity (\ref{mainids4}) as an equation with respect to $\hat w_4$. The solution $\hat w_4=w_4$ leads to the point $W$ (see (\ref{expr234})). Further, the quadratic equation  
\begin{equation}
\begin{array}{l}
(i-1)w_3^2\,\hat w_4^2+ \Bigl((1-2i)w_3+(i-1)w_3^2w_4\Bigr)\hat w_4+\\
\vspace{-0.3cm}\\
\hspace{5cm}\Bigl(i+(1-2i)w_3w_4+(i-1)w_3^2w_4^2\Bigr)=0\label{qudreq}
\end{array}
\end{equation}
has the following solutions:
\begin{equation}
\hat w_4=\frac{2i-1+(1-i)w_3w_4+\sqrt{6iw_3^2w_4^2-(2+6i)w_3w_4+1}}{(2i-2)w_3}.\label{solquadr}
\end{equation}
Our goal now is to choose $\varepsilon$ so that for $W\in U_{\varepsilon}$ neither of the points $\hat W$ defined by solutions (\ref{solquadr}) lies in $U_{\varepsilon}$. We write $w_2=\bar w_1+\mu$, $w_4=\bar w_3+\eta$, $\hat w_4=\bar w_3+\hat\eta$ and show that $\varepsilon$ can be taken so that if $|\mu|<\varepsilon$, $|\eta|<\varepsilon$, then $|\hat\eta|\ge\varepsilon$. 

Formula (\ref{solquadr}) implies
\begin{equation}
\begin{array}{l}
-8|w_3|^2+4+i(24|w_3|^4-24|w_3|^2+4)+\\
\vspace{-0.1cm}\\
\hspace{1cm}\Bigl[24i\eta|w_3|^2w_3+8i\eta^2w_3^2-(4+12i)\eta w_3\Bigr]=\\
\vspace{-0.1cm}\\
\hspace{3cm}-24i\hat\eta|w_3|^2w_3-8i(\hat\eta^2+\eta\hat\eta)w_3^2+(4+12i)\hat\eta w_3.
\end{array}\label{mainids1}
\end{equation}
Observe that for any $w_3$ one has
\begin{equation}
\Bigl|-8|w_3|^2+4+i(24|w_3|^4-24|w_3|^2+4)\Bigr|\ge\frac{4}{3}.\label{cond1}
\end{equation}
Next, since $\varepsilon<1$, formula (\ref{formmm18}) yields
$$
\left|\bar w_1+\frac{\mu}{2}\right|^2+\left|\bar w_3+\frac{\eta}{2}\right|^2=1+\frac{|\mu|^2}{4}+\frac{|\eta|^2}{4}<1+\frac{\varepsilon^2}{2}<\frac{3}{2},
$$
which implies 
\begin{equation}
|w_3|<2.\label{est4}
\end{equation}
Using (\ref{est4}), one can estimate the terms in square brackets in the left-hand side of (\ref{mainids1}) as follows:
\begin{equation}
\Bigl|24i\eta|w_3|^2w_3+8i\eta^2w_3^2-(4+12i)\eta w_3\Bigr|< 32\varepsilon^2+224\varepsilon.\label{cond2}
\end{equation}
Similarly, using (\ref{est4}) for the right-hand side of (\ref{mainids1}) we have
\begin{equation}
\Bigl|-24i\hat\eta|w_3|^2w_3-8i(\hat\eta^2+\eta\hat\eta)w_3^2+(4+12i)\hat\eta w_3\Bigr|<32|\hat\eta|^2+256|\hat\eta|.\label{cond3}
\end{equation}

It follows from formulas (\ref{mainids1}), (\ref{cond1}), (\ref{cond2}), (\ref{cond3}) that
\begin{equation}
32|\hat\eta|^2+256|\hat\eta|>\frac{4}{3}-(32\varepsilon^2+224\varepsilon).\label{ineq}
\end{equation}
Thus, in order to finalize the proof of the theorem, we need to choose $\varepsilon$ so that inequality (\ref{ineq}) implies $|\hat\eta|\ge\varepsilon$. For example, let $\varepsilon$ be such that
$$
32\varepsilon^2+224\varepsilon<\frac{1}{3}.
$$
For instance, $\varepsilon=\sqrt{2}\cdot 10^{-3}$ satisfies this condition. Then from (\ref{ineq}) one has
$$
32|\hat\eta|^2+256|\hat\eta|>1,
$$
which implies $|\hat\eta|>0.003>\varepsilon$ as required. The proof of the theorem is complete. \qed

\begin{remark}\label{remimpr}\rm By experimenting with inequality (\ref{ineq}) one can slightly improve the value of $\varepsilon$. However, the improved value is still of order $10^{-3}$, thus the lower bound for $t_0-1$ it leads to is still of order $10^{-6}$.
\end{remark}

We finish the paper by making several observations regarding solutions (\ref{solquadr}) of equation (\ref{qudreq}) concentrating on the case when $W\in M_t^3$ for some $t<\sqrt{(2+\sqrt{2})/3}$.
\vspace{0.3cm}

(i) It follows from (\ref{phieq}) that $w_4$ is a solution of (\ref{qudreq}) if and only if $J_{\tilde F}$ vanishes at the point $W$. Therefore, by Proposition \ref{embed}, neither of the values in the right-hand side of (\ref{solquadr}) is equal to $w_4$ if $W\in M_t^3$ with $t<\sqrt{(2+\sqrt{2})/3}$. 
\vspace{0.3cm}

(ii) Arguing as in the proof of Proposition \ref{embed}, one can see that the polynomial $6iw_3^2w_4^2-(2+6i)w_3w_4+1$ vanishes at some point of $M_t^3$ if and only if $t\ge 2/\sqrt{3}$. Since $2/\sqrt{3}>\sqrt{(2+\sqrt{2})/3}$, it follows that the right-hand side of (\ref{solquadr}) defines two distinct points $\hat W_1$, $\hat W_2$ if $W\in M_t^3$ with $t<\sqrt{(2+\sqrt{2})/3}$. Combining this fact with (i), we see that for such $W$ the fiber of $\tilde F$ over $\tilde F(W)$ consists of exactly three points.  
\vspace{0.3cm}

(iii) Using formula (\ref{solquadr}) one can determine, in principle, all values\linebreak $t<\sqrt{(2+\sqrt{2})/3}$ such that neither of $\hat W_1$, $\hat W_2$ lies in $M_t^3$ provided $W\in M_t^3$, which is equivalent to the injectivity of the map $\tilde F$ on  $M_t^3$. However, the computations required for this analysis are rather complicated, and we did not carry them out in full generality. These computations significantly simplify if, for instance, $w_4=0$. In this case formula (\ref{solquadr}) yields the solutions
\begin{equation}
\hat w_4=\frac{1}{w_3}, \quad\hat w_4=\frac{1-i}{2w_3}.\label{twoptsspecial}
\end{equation}        
Also, we have
$$
|w_1|^2+\frac{1}{|w_1|^2}+|w_3|^2=2t,
$$
which implies $|w_3|<|w_1|$ since otherwise $t\ge\sqrt{2}>\sqrt{(2+\sqrt{2})/3}$. Therefore, for each of the two points $\hat W_1$, $\hat W_2$ given by (\ref{twoptsspecial}) one has
$$
|w_1|^2+|\hat w_2|^2+|w_3|^2+|\hat w_4|^2>2t,
$$
hence neither of these points lies in $M_t^3$. We also arrived at the same conclusion when analyzing several other special cases. Thus, our calculations lead to the following conjecture.

\begin{conjecture}\label{conj} \sl The restriction of $\tilde F$ to $M_t^3$ is an embedding for all parameter values $1<t<\sqrt{(2+\sqrt{2})/3}$.
\end{conjecture}

(iv) We note that the non-injectivity of $\tilde F$ on $M_t^3$ is easy to see, for example, if $t\ge\sqrt{2}$. Indeed, let $u\ne 0$ be a real number satisfying
$$
2u^2+\frac{1}{u^2}=2t, 
$$
and consider the following three distinct points in $Q^3$:
\begin{equation}
\hspace{-0.15cm}W_u:=\left(u,\frac{1}{u},u,0\right), W_u':=\left(u,0,u,\frac{1}{u}\right), W_u'':=\left(u,\frac{1+i}{2u},u,\frac{1-i}{2u}\right).\label{threepoints}
\end{equation}
Then $W_u, W_u', W_u''\in M_t^3$ and $\tilde F(W_u)=\tilde F(W_u')=\tilde F(W_u'')=(u,u,0)$. Since every fiber of $\tilde F$ contains at most three points, $W_u,W_u', W_u''$ form the complete fiber of $\tilde F$ over the point $(u,u,0)$.
\vspace{0.3cm}

(v) In \cite{AR} the authors in fact introduced not just the map $f$ (see (\ref{mapf})) but a class of maps of the form
$$
g:\CC^2\ra\CC^3,\quad g(z,w):=(z,w,P(z,\bar z,w,\bar w)).
$$
Here $P$ is a harmonic polynomial given by
$$
P=\left(\bar z\frac{\partial}{\partial w}-\bar w\frac{\partial}{\partial z}\right)\left(\sum_{j=1}^m\frac{1}{p_j(q_j+1)}Q_j\right),
$$
where $Q_j$ is a homogeneous harmonic complex-valued polynomial on $\CC^2$ of total degree $p_j\ge 1$ in $z,w$ and total degree $q_j$ in $\bar z$, $\bar w$, such that the sum $Q:=Q_1+\dots+Q_m$ does not vanish on $S^3$. Observe that although the third component of the map $f$ is not harmonic, it is obtained (up to a multiple) from the harmonic polynomial $w\bar z\bar w^2-z\bar z^2\bar w+i\bar z\bar w$ by replacing the term $i\bar z\bar w$ with $i\bar z\bar w(|z|^2+|w|^2)$ (these two expressions coincide on $S^3$). It is convenient to take $Q$ to be a polynomial in $|z|^2$, $|w|^2$ (as was done in \cite{AR}), but in this case $P$ is divisible by $\bar z\bar w$, which implies that the holomorphic extension $G$ of the push-forward of $g$ to $\RR^4$ is not injective on $M_t^3$ with $t\ge\sqrt{2}$. Indeed, writing $G$ in the coordinates $w_j$ defined in (\ref{coordw}), for the points $W_u$, $W_u'$ introduced in (\ref{threepoints}) one has $G(W_u)=G(W_u')=(u,u,0)$. Thus, one cannot obtain the embeddability of $M_t^3$ in $\CC^3$ for $t\ge\sqrt{2}$ by utilizing any of the maps introduced in \cite{AR} with $Q$ being a function of $|z|^2$, $|w|^2$ alone.

{\obeylines
Department of Mathematics
The Australian National University
Canberra, ACT 0200
Australia
e-mail: alexander.isaev@anu.edu.au
}

\end{document}